\documentclass[10pt,letterpaper]{article}

\usepackage{amsmath,amssymb}

\usepackage{changepage}

\usepackage[utf8x]{inputenc}

\usepackage{textcomp,marvosym}

\usepackage{cite}

\usepackage{nameref,hyperref}

\usepackage[right]{lineno}

\usepackage{microtype}
\DisableLigatures[f]{encoding = *, family = * }

\usepackage[table]{xcolor}

\usepackage{array}

\newcolumntype{+}{!{\vrule width 2pt}}

\newlength\savedwidth

\usepackage{float}

\usepackage[aboveskip=1pt,labelfont=bf,labelsep=period,justification=raggedright,singlelinecheck=off]{caption}

\makeatletter
\renewcommand{\@biblabel}[1]{\quad#1.}
\makeatother

\date{}

\usepackage{lastpage,fancyhdr,graphicx}

\begin{document}
\vspace*{0.2in}

\begin{flushleft}
{\Large
\textbf\newline{Admission predictors for  success in a mathematics graduate program} 
}
\newline
\\
Timmy Ma\textsuperscript{1\Yinyang},
Karen E. Wood\textsuperscript{1\Yinyang},
Di Xu\textsuperscript{2},  
Patrick Guidotti, \textsuperscript{1}, 
Alessandra Pantano, \textsuperscript{1}, 
Natalia L. Komarova\textsuperscript{1*}
\\
\bigskip
\textbf{1} Department of Mathematics, University of California, Irvine, Irvine, CA, 92617
\\
\textbf{2} School of Education, University of California, Irvine, Irvine, CA, 92617

\bigskip

%
%
\Yinyang These authors contributed equally to this work.





* komarova@uci.edu

\end{flushleft}
\section*{Abstract}
There are many factors that can influence the outcome for students in a mathematics PhD program: bachelor's GPA (BGPA), bachelor's major, GRE scores, gender, Under-Represented Minority (URM) status, institution tier, etc. Are these variables equally important predictors of a student's likelihood of succeeding in a math PhD program? In this paper, we present and analyze admission data of students from different groups entering a math PhD program at a southern California university. We observe that some factors correlate with success in the PhD program (defined as obtaining a PhD degree within a time-limit). According to our analysis, GRE scores correlate with success, but interestingly, the verbal part of the GRE score has a higher predictive power compared to the quantitative part. Further, we observe that undergraduate student GPA does not correlate with success (there is even a slight negative slope in the relationship between GPA and the probability of success). This counterintuitive observation is explained once undergraduate institutions are separated by tiers: students from ``higher tiers" have undergone a more rigorous training program; they on average have a slightly lower GPA but run a slightly higher probability to succeed. Finally, a gender gap is observed in the probability to succeed with female students having a lower probability to finish with a PhD despite the same undergraduate performance, compared to males. This gap is reversed if we only consider foreign graduate students. It is our hope that this study will encourage other universities to perform similar analyses, in order to design better admission and retention strategies for Math PhD programs.


\section*{Introduction}
A higher education degree, in particular in STEM disciplines, is a desirable goal for many millenials; for example, it has been shown that there is a direct correlation between increased earnings and quality of life on the one hand,  and the pursuit of a graduate degree  on the other \cite{baum2010education, baum2013education}. 
However, students who do not ``match" well with the program  are less likely to graduate and to benefit from attending graduate school. The process of pursuing/granting a PhD degree involves significant costs for both the individual students and the institutions;  as a result, there is a need for universities to identify applicants who are most promising in their PhD  program. In addition, understanding the predictors of success could potentially help institutions with the early identification of students who are encountering more challenges, and design ways to better support them. Female and underrepresented minorities (URM) comprise a particularly vulnerable population in graduate programs across the STEM fields. An effort is underway to understand the causes of  underrepresentation and implement measures to mitigate it. This study seeks to identify factors that may or may not  predict success in a mathematics PhD program. \\

In this study, we gathered and analyzed data from a graduate program in Mathematics of about 100 students at a public Tier I university in southern California. We were particularly interested in detecting any discrepancy between the success rates of students of different groups which could not be explained by differences in their undergraduate record. \\

In this institution, PhD students are admitted at a rate of about 15-25 per year from an applicant pool of about 200-300. Of the admitted students, two thirds are domestic or permanent residents and the rest are international students, mainly from China. Admitted students are guaranteed six years of financial support while in good academic standing.
The average degree completion time is 5-6 years,  while the maximum time permitted is seven years. The program is fairly standard and consists of preparatory and specialization courses, two qualifying examinations that can be chosen from a set of three (algebra, real, and complex analysis), and a thesis on a research topic. The overall success rate in the program is about 55\% as computed over the most recent period of six years (admissions from 2005 to 2010) for which students have reached maximum allowed time to degree. For comparison, the national average success rate in graduate programs in Mathematics as reported by the Council of Graduate Schools Ph.D Completion project is 47.5\% \cite{PhDComplete}. 
While the  sample under consideration is limited, the results obtained can potentially inform comparable graduate programs and encourage them to perform similar analyses. 

\section*{Data}

\begin{figure}[H]
\begin{adjustwidth}{-1.25in}{0in}
\centering
\caption{Data demographics: (left) Plot of the population by Citizenship and colored by Gender. (right) Plot of the population by Ethnicity and colored by Gender.\label{PopCitGen}}
\includegraphics[scale=0.8]{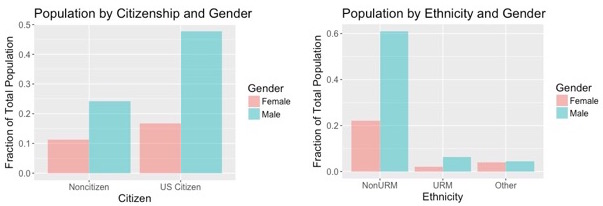}
\end{adjustwidth}
\end{figure}

\begin{figure}[H]
\begin{adjustwidth}{-1.25in}{0in}
\centering
\caption{Plot of the student population  by the three cohorts of collected data and colored by Gender. \label{PopRangeGen}}
\includegraphics[scale=0.8]{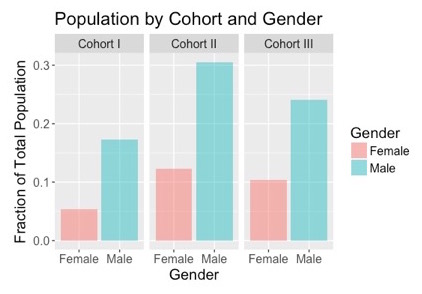}
\end{adjustwidth}
\end{figure}

\begin{figure}[H]
\begin{adjustwidth}{-1.25in}{0in}
\centering
\vskip.5cm
\caption{Histogram comparing undergraduate performance of the male and female students in the program, as measured by BGPA, MGPA, GREV, and GREQ. \label{All_data}}
\includegraphics[scale=0.8]{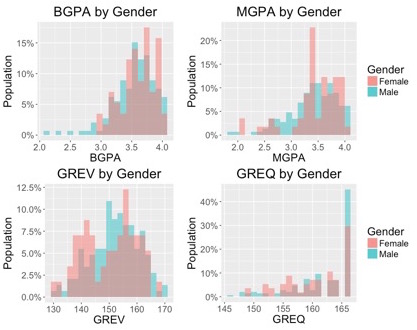}
\end{adjustwidth}
\end{figure}

\begin{figure}[H]
\begin{adjustwidth}{-1.25in}{0in}
\centering
\vskip.5cm
\caption{Histogram  comparing undergraduate performance of the URM students in the program, as measured by BGPA, MGPGREV, and GREQ. \label{All_data_2}}
\includegraphics[scale=0.8]{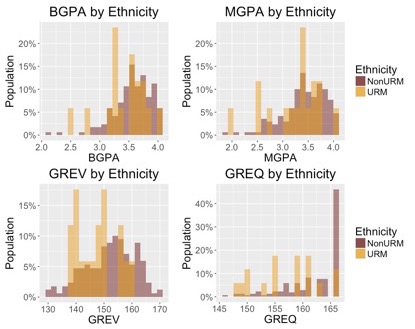}
\end{adjustwidth}
\end{figure}

\begin{figure}[H]
\begin{adjustwidth}{-1.25in}{0in}
\centering
\vskip.5cm
\caption{Histogram comparing undergraduate performance of students by Citizenship in the program: BGPA, MGPA, GREV, and GREQ. \label{All_data_3}}
\includegraphics[scale=0.8]{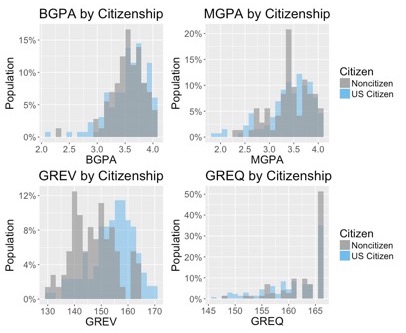}
\end{adjustwidth}
\end{figure}

We have collected student data (described below), from a PhD program in mathematics, at a public institution in Southern California. Our database contains 203 students from  nine consecutive admissions years. Each of these students were  followed  for 7 years (the maximum time allowed for graduation  by the PhD program). For simplicity, in this study we have split our student sample into three cohorts based on the year of admission (Year 1-3, 4-6, 7-9).\footnote{We did not include students  who were admitted to the Masters degree program, as they were not pursuing a PhD.}  Table \ref{tab_all} contains information about the  demographics of the students  (gender, citizenship, URM status and undergraduate major), as well as information about their undergraduate performance in mathematics (as measured by the overall and math Bachelor's GPA, and the verbal and quantitative GRE scores).  Note that 
 the  mathematics GRE subject test was not required for admission to this PhD program during the time-span under consideration. Fig \ref{PopCitGen}-\ref{All_data_3}  break down this information into finer categories.

\begin{table}[H]
\begin{adjustwidth}{-1.25in}{0in}
\caption{\label{tab_all}Demographics of the data and mean undergraduate performance by category: bachelor's GPA (BGPA), Math undergraduate GPA (MGPA), GRE verbal (GREV), and GRE quantitative (GREQ).}
\begin{center}
\begin{tabular}{|l|c|c||c|c|c|c|}
\hline
Category 		& \# students & \% & Mean BGPA 	&  Mean MGPA	& Mean GREV & Mean GREQ \\
\hline
Male 			& 146 	&	71.9 	& 3.55 	& 3.41 	& 152.64 	& 161.81\\
Female			& 57 	&	28.1	& 3.61	& 3.44 	& 150.72 	& 160.54 	 \\
\hline
US Citizens 	& 131 	& 	64.5	& 3.57 	& 3.41 	& 154.51 	& 160.66\\
Non US 			& 72 	&	35.5	& 3.56 	& 3.43 	& 147.71 	& 162.91 \\
\hline
Non-URM			& 169 	& 	83.3 	& 3.57 	& 3.43 	& 152.63	& 162.09\\
URM				& 17 	& 	8.4 	& 3.43 	& 3.24	& 147.94	& 157.06\\
Other/Unknown 	& 17	& 	8.4 	& 3.60	& 3.44 	& 151.00 	& 159.59\\
\hline
Math Major 		& 178 	& 	87.7 	& 3.57	& 3.42 	& 151.75 	& 161.48\\
Non Math Major 	& 13 	& 	6.4		& 3.45	& 3.42	& 156.08 	& 162.53	\\
NAs (Unknown)	& 12 	& 	5.9 	& 3.61	& 3.34 	& 152.92 	& 159.92	\\
\hline
\end{tabular}
\end{center}
\end{adjustwidth}
\end{table}


The last four columns of Table \ref{tab_all} present numerical information about the students' undergraduate performance, split by category. The mean value of bachelor's GPA (BGPA), Math undergraduate GPA (MGPA), GRE verbal (GREV), and GRE quantitative (GREQ) across each category is provided.
Not all differences in the mean scores are significant, as explained below. There are three sets of comparisons, based on the categories. 
\begin{itemize}
\item {\bf Gender.} In Fig \ref{All_data}, undergraduate GPA and  GRE scores are split by gender. We can see that males and females perform comparably in terms of both GPA and GRE. The t-test reveals no statistically significant difference between  male and female mean scores in BGPA ($p=0.176$), MGPA ($p=0.585$),  for GREV ($p=0.180$), or for GREQ ($p=0.115$).

\item {\bf URM status.} In Fig\ref{All_data_2},  we split the GPA and GRE data by the URM status. The t-test shows that there is no statistically significant difference in performance between URM and non-URM students  for BGPA ($p=0.1571$) or for MGPA ($p=0.1638$). However, URM students have a slightly lower GRE mean ($p=0.0184$ for GREV, and  $p=0.0026$ for GREQ).

\item {\bf Nationality.} In Fig \ref{All_data_3}, the GPA and GRE data are split by  citizenship status. We observe that the undergraduate GPA is similar between US and non-US citizens ($p=0.864$ for BGPA and $p=0.6995$ for MGPA). As for GRE, international students perform differently fom US students, reporting a lower mean in GREV ($p<0.0001$) and a higher mean in GREQ ($p=0.001$). 
\end{itemize}

Information on the students' undergraduate institution was categorized as follows. PhD granting US institutions were divided in tiers I/II/III, according to the AMS ranking for PhD granting institutions (Private and Public); institutions whose highest degree is Bachelors and Masters were  included in categories B and M, respectively; Chinese institutions were  grouped according to their ranking and split into three categories (top twenty schools,  schools ranked from 21-100, and schools ranked below the top 100). Category ``Other" contained students from non-US and non-Chinese institutions. These data are summarized in Table \ref{tab:tiers}. 

\begin{table}[H]
\begin{adjustwidth}{-1.25in}{0in}
\caption{\label{tab:tiers}Table of undergraduate college tier levels.}
\begin{center}
\begin{tabular}
{|m{.15\textwidth}|m{.15\textwidth}|m{.1\textwidth}|m{.15\textwidth}|m{.1\textwidth}|m{.15\textwidth}|}
\hline
\multicolumn{2}{|m{.3\textwidth}}{AMS Rankings \newline (PhD-granting US Schools)} & \multicolumn{2}{|m{.25\textwidth}}{Lettered Groups} & \multicolumn{2}{|m{.25\textwidth}|}{Numerical Tiers \newline (Chinese Institutions)}  \\
\hline 
Ranking & \# students & Ranking & \# students & Ranking & \# students\\
\hline
\textbf{I - Public} \newline \textbf{I - Private} & 40 \newline 5 & \textbf{B} & 19 & \textbf{1-20} & 34 \\
\textbf{II} & 48 & \textbf{M} & 26 & \textbf{21-100} & 6 \\
\textbf{III} & 2 & \textbf{Other} & 20 & \textbf{100+} & 3\\
\hline
\end{tabular}
\end{center}
\end{adjustwidth}
\end{table}

\section*{Analysis and Results}
The purpose of this analysis is to gain a perspective of what factors, if any, predict  the success of a student in a math PhD program. We use a general logistic regression model for predicting the qualitative response $Y$ on the basis of several predictor variables $X_i$:

\begin{eqnarray}
\label{regression_eq}
Y = \beta_0 + \beta_1 X_1 + \beta_2 X_2 + \dots + \beta_n X_n + \epsilon,
\end{eqnarray}
where $\epsilon$ is a mean-zero random error term, and $\beta_i$ quantifies the association between variable $X_i$ and response $Y$. We interpret $\beta_i$ as the average effect on $Y$ of a one unit increase in $X_i$, holding all other predictors fixed. A positive (negative) $\beta_i$ suggests a correlation (an anti-correlation) between predictor $X_i$ and response $Y$. 

In this study, $Y$ is ``Success in the PhD program", defined as PhD completion within 7 years. Any other outcome, such as leaving the PhD program after failing qualifying exams, staying longer than 7 years (there were only two students in this scenario), or leaving the PhD program after receiving a MS degree amounts to failure (lack of success) in our study. \footnote{For the purpose of this study we chose to focus exclusively on success in obtaining a PhD degree, there are of course other valid metrics for success.} 

\begin{figure}[H]
\begin{adjustwidth}{-1.25in}{0in}
\caption{Plot of regression for BGPA,  in the absence of Tier information (left) and when Tiers are taken into account (right). The x-axis represents BGPA, the y-axis is prediction of success with the regression model.\label{BGPA}}
\centering
\includegraphics[scale=0.8]{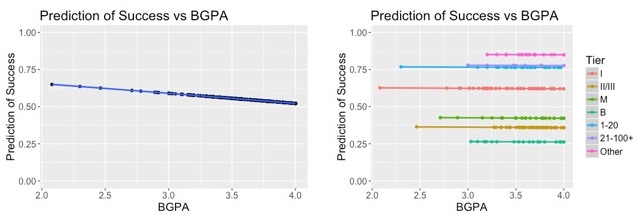}
\end{adjustwidth}
\end{figure}

The first interesting finding comes from computing the regression on the BGPA \emph{only}. The result is counterintuitive, as Bachelor's GPA turns out to have a negative correlation  to the success of students in the PhD program (coefficient $=-0.27$, $p=0.538$), see Fig \ref{BGPA} (left). A possible explanation could be that the same grade point average obtained in two programs with different levels of academic rigor may reflect a different level of students' preparedness. Therefore, we performed the regression again, this time taking into consideration both the BGPA and the type (Tier) of undergraduate institution  (coefficient $=-0.014$, $p=0.977$). The results are presented in Fig \ref{BGPA} (right).\footnote{In decreasing order of likelihood of success, Fig \ref{BGPA} (right) shows tiers:  Other,  1-20, 21-100/100+, I, M, II/III and B.} We observe  that  differences in BGPA are overridden by the Tier. Moreover, students obtaining their bachelors' from foreign institutions are predicted to do better than students graduating from domestic institutions (coefficient $=0.073$, $p=0.935$). Within US institutions, schools of Tier I perform the best and bachelor granting institutions the worst (coefficient $=-2.208$, $p<0.001$).

\begin{figure}[H]
\begin{adjustwidth}{-1.25in}{0in}
\caption{Plot of regression for GREV (left) and GREQ (right) colored by Tier. The x-axis represents the  GRE score (GRE verbal on the left and GRE quantitative on the right), the y-axis is prediction of success with the regression model.\label{GREV}}
\centering
\includegraphics[scale=0.8]{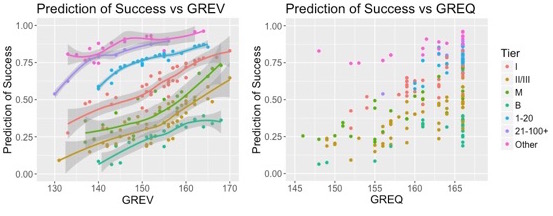}
\end{adjustwidth}
\end{figure}

Next, we perform regression on only the GREV with Tier, and only on 
 GREQ with Tier. The results are plotted in Fig \ref{GREV}.
 We observe  a clear positive correlation between the score on the GRE verbal test and prediction of success in the PhD program (coefficient $=0.051$, $p=0.019$).  Moreover, certain tiers of schools are performing better than others. In regards to GRE quantitative test, we continue to see a positive correlation between GREQ score and prediction of success  (coefficient $=0.065$ and  $p=0.07$),  but this correlation is only  marginally significant, possibly due to the small variation in the GREQ scores among math PhD students (all scores tend to be near the maximum possible score of 166, with about a half within 2 points, and three quarters within 7 points of the maximum).
 
To perform a systematic analysis that includes all variables of interest, we studied regression with all the  observable characteristics of students, see Tables \ref{tab3} and \ref{tab4}.  The $R^2$ value increases from 0.0021 to 0.2307 as we introduce more categories, which means that a better fit is achieved for the data. In particular, the correlation between GREQ and success becomes stronger once controlling for other baseline characteristics. \\

\begin{table}[H]
\begin{adjustwidth}{-1.25in}{0in}
\caption{\label{tab3}Table of coefficients for regression. Column ME (SE) are the values of marginal effects and the standard error. The significance codes are: $0$ `***' $0.001$ `**' $0.01$ * $0.05$ `.' $0.1$ ` ' 1.}
\begin{center}
\begin{tabular}{|l|c|c|c|c|c|c|c|c|c|c|}
\hline
Category	&	[1]		&	[2]		& [3]		& [4]		& [5]		& [6] 	& ME (SE)	\\
\hline
Intercept	&	0.03509	& 0.6082 .	& 1.4364	& -5.18925	&-13.9925**	&	-14.0524** & \\
\hline
Gender (Male) & 0.24061 & 0.3187	& 0.3039	& 0.2418	& 	0.1633	& 	0.1607	& 0.039 (0.08)\\
\hline
Citizenship (US) &  	& -0.9520** & -0.9498***& -1.298***& -1.04400**	& 	-1.0479** & -0.247 (0.08)\\
\hline
BGPA 		& 			& 			& -0.2299 	& -0.2914	& 	-0.3726	& 	-0.5048	& -0.124 (0.144)\\
\hline
GREV 		& 			& 			&	  		& 0.0468*	& 	0.0311	& 	0.03179 & 0.0078 (0.005)	\\
\hline
GREQ 		& 			& 			&  			& 			&   0.0705*	& 	0.0697* & 0.017 (0.008)	\\
\hline
MGPA		& 			& 			& 			& 			& 			& 	0.1637	& 0.040 (0.104)\\
\hline \hline
$R^2$ 		& 0.002111 	& 0.03709 	& 0.03800	& 0.0606	& 	0.0784	& 	0.0790& 	\\
\hline
\end{tabular}

\end{center}
\end{adjustwidth}
\end{table}

\begin{table}[H]
\begin{adjustwidth}{-1.25in}{0in}
\begin{center}
\caption{\label{tab4}Table of coefficients for regression, continued. Column ME (SE) are the values of marginal effects and the standard error. The significance codes are: $0$ `***' $0.001$ `**' $0.01$ * $0.05$ `.' $0.1$ ` ' 1.}
\begin{tabular}{|l|c|c|c|c|c|}
\hline
Category		& [7]		& [8]		& [9]		& [10]	& ME (SE)\\
\hline
Intercept		& -14.571**	& -11.075*		& -33.036	& -33.031	&	\\
\hline
Gender (Male) 	& 0.1441	& 0.0878		& 0.1312	& 0.1273	& 0.030 (0.374)\\
\hline
Citizenship (US) &-1.015**	&  -0.9894*		& 17.039	& 17.139	&0.997 (1.327)\\
\hline
BGPA 			&-0.5594	& -0.6269		& -0.6027	& -0.6172	&-0.148 (1.672)\\
\hline
GREV 			& 0.0310	& 0.02702		& 0.0382	& 0.0373	&0.0089 (0.101)\\
\hline
GREQ 			& 0.0731*	& 0.5891 .		& 0.0696 . 	& 0.0715 .	& 0.017 (0193)\\
\hline
MGPA			& 0.1906	& 0.1572		& 0.2506	& 0.2531	&0.061 (0.692)\\
\hline
Major (Math) 	& 0.1527	& 0.1051		& 0.8411 	& 0.8204	& 0.178 (3.27)\\
Major (Unknown) & 0.6201	& 0.6886		& 0.7733	& 0.7427	& 0.183 (0.517)\\
\hline
Ethnicity (URM)	& 			& -1.039 .		& -1.133	& -1.052	& -0.217 (4.517)\\
Ethnicity (Other)	& 		& -1.066 .		& -1.532*	& -1.422* 	& -0.274 (6.298)\\
\hline
Tier (II/III)	& 			& 				& -0.7558	& -0.777	&-0.174 (2.739)\\
Tier (M)		& 			& 				& -0.4210	& -0.417 	&-0.096 (1.454)\\
Tier (B)		& 			& 				& -1.8558**	& -1.860 ** & -0.3302 (8.128)\\
Tier (1-20)		& 			& 				& 17.396	& 17.503 	&0.9662 (8.480)\\
Tier (21-100)	& 			& 				& 17.745	& 17.94 	& 0.771 (19.160)\\
Tier (Other)	& 			& 				& 19.569	& 19.62		& 0.9133 (13.862)\\
\hline
Year II 		&			&				&			& -0.3679	&-0.087 (1.041)\\
Year III		&			&				&			& -0.1652	&-0.039 (0.481)\\
\hline \hline
$R^2$ 			& 0.0809	& 	0.1032		& 	0.2281	& 	0.2307 & \\
\hline
\end{tabular}
\end{center}
\end{adjustwidth}
\end{table}

\begin{figure}[H]
\begin{adjustwidth}{-1.25in}{0in}
\centering
\caption{Plot of prediction for Tiers from regression of Table \ref{tab4}, column [10]. The x-axis represents the Tier level, y-axis is prediction of success with the regression model. The top and bottom of the box indicate the 3rd and 1st quartile of the data, respectively. The line in the middle of the box is the median of the data. The ``whiskers" are the maximum and minimum values without outliers, and the dots represent outliers. Note that the latter are values falling outside 1.5 times the interquartile range above the 3rd and below the 1st quartile range. \label{Tier}}
\includegraphics[scale=0.8]{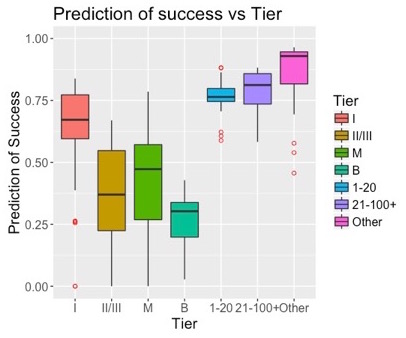}
\end{adjustwidth}
\end{figure}

Fig \ref{Tier} shows the predicted success rate for students from different tiers using the regression of  Table \ref{tab4}, column [10]. We observe some descriptive differences among tiers. For example, students from tier I schools are predicted to perform better than students from tier II/III schools or students from master degree granting schools. International students are predicted to perform the best of all, and students from bachelor's granting institutions the worst. None of the differences however reach statistical significance except for Bachelor's degree only granting institution.\\

\begin{figure}[H]
\begin{adjustwidth}{-1.25in}{0in}
\centering
\caption{Plot of prediction for Gender:  all data from regression of  Table \ref{tab4}, column [10] (left); split by Tier  (right). \label{Gender}}
\includegraphics[scale=0.8]{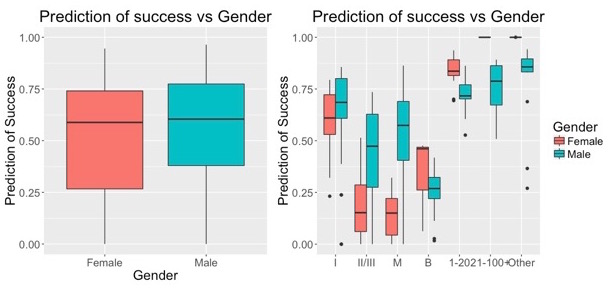}
\end{adjustwidth}
\end{figure}

Fig \ref{Gender} (left) shows the success rate for students by gender using the regression of  Table \ref{tab4}, column [10]. We notice a slight gender gap in the data, with males predicted to succeed better than females (coefficient $=0.127$, $p>0.1$). To test whether this gender gap could be due to some interactional effects between gender and type of undergraduate institution, we added interaction terms between Gender and Tier. The prediction is plotted in Fig \ref{Gender} (right). We can see that the gender gap is different in different tiers. None of the interactions however reach statistical significance, possibly due to the small sample size. \\

\begin{figure}[H]
\begin{adjustwidth}{-1.25in}{0in}
\centering
\caption{Plot of prediction for Citizenship from regression of column [10]. x-axis is the Citizenship status, y-axis is prediction of success with the regression model. (Left) All data, (right) by gender. \label{Citizen}}
\includegraphics[scale=0.8]{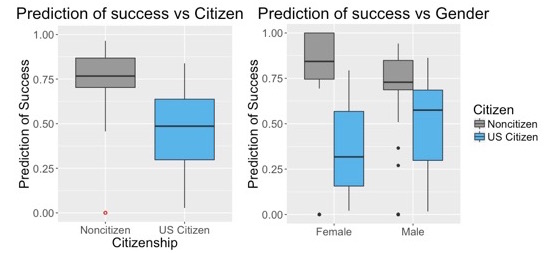}
\end{adjustwidth}
\end{figure}

Fig \ref{Citizen} (left) shows the success rate for students by citizenship status using the regression of Table \ref{tab4}, column [10]. We can see a gap where Non-US students are predicted to succeed at a higher rate than US students (coefficient $=17.139$, $p>0.1$). We should note that citizenship is highly correlated with tier (e.g. we separated tier dummies to indicate institutions in Europe and China). As such, the raw correlation between citizenship and success (Fig \ref{Citizen}) is absorbed in the statistical analysis that controls for other variables. Fig \ref{Citizen} (right), for example, shows  the correlation between Gender and Citizenship. We can see that Non-US female students are more likely to succeed than (all) males, but  US male students  are more likely to succeed than US female students. We observe a similar pattern when we look at tiers in Fig \ref{Gender} (right).

\begin{figure}[H]
\begin{adjustwidth}{-1.25in}{0in}
\centering
\caption{Plot of prediction for Ethnicity from regression of column [10]. x-axis are the Ethnicities, y-axis is prediction of success with the regression model. \label{Eth}}
\includegraphics[scale=0.8]{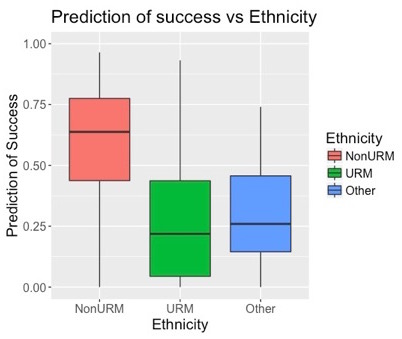}
\end{adjustwidth}
\end{figure}

Fig \ref{Eth} shows the success rate for students of different ethnicities and demonstrates the existence of a clear racial gap. When performing regression  using column [8] in Table \ref{tab4}, this result appears to be only marginally significant ($p = 0.086$). However, the lack of statistical significance is more likely due to a lack of power to detect it precisely - because of the small sample size, rather than to a lack of relationship between race and success. The same reasoning  applies to other variables that clearly demonstrate a difference in coefficients, but do not reach statistical significance. 

\section*{Discussion}

In this study we analyzed the success rate of  of PhD students admitted to a math graduate program of a public Tier I institution in southern California in 9 consecutive years (for a total 203 students). The following patterns have been observed. 

\begin{itemize}

\item There is a positive correlation between GREV/GREQ scores and probability of success in the PhD program. The correlation is stronger for the GRE verbal than for GRE quantitative, due to the small variation of GREQ scores among math PhD students. 

\item Despite the fact that the undergraduate performance of male and female students is comparable (in terms of GRE and BGPA), a gender gap exist in their success rate in the PhD program. Similar results hold for URM/non-URM students. 

\item  Citizenship status (foreign/US citizen) is a predictor of success. Gender gap is observed among all students, and persists among  US students: male US students are predicted to do better than female US students. Female foreign students on the other hand are predicted to do better than their male counterparts. 

\item Counterintuitively, undergraduate GPA negatively correlates with success in the graduate program.
This phenomenon is explained by splitting the undergraduate institutions into tiers. Students that come from ``higher tiers" (which possibly translates into a more rigorous undergraduate program and a lower GPA) tend to have a higher chance to succeed. 

\item Gender gap persists within each tier of undergraduate institutions, but it is different in different tiers. 

\end{itemize}

Although we did not speak to the specifics of why students leave the PhD program (due to time restraints, passing of qualifying exams, making adequate progress in research, or finding a high paying job), this is an important aspect that should be examined once the data become available. 

This study adds to a growing effort to understand what correlates to the success in PhD programs. In Ref. \cite{moneta2017limitations}, GRE scores, undergraduate GPA, and other factors about the undergraduate students entering a biomedical PhD program were studied; it was found neither GREV nor GREQ scores  predicted success, measured as graduating with a PhD. Undergraduate GPA on the other hand was shown to significantly predict graduation with a Ph.D. Ref. \cite{pacheco2015beyond} proposed a composite score that appeared to be a better predictor than the GRE in the context of a biomedical PhD program. Differences in these and our findings are interesting and may be attributed to the differences in the fields of study (biomedicine vs mathematics), but clearly more research must be performed in order to understand these issues better.

\section*{Acknowledgments}

\nolinenumbers

%
%
%

\end{document}